\documentclass[a4paper]{amsart}

\usepackage[british]{babel}

\usepackage{mlmodern}
\DeclareFontFamily{OMX}{mlmex}{}
\DeclareFontShape{OMX}{mlmex}{m}{n}{<->mlmex10}{}
\usepackage[T1]{fontenc}
\usepackage[babel,stretch=11,shrink=11,verbose=errors,selected]{microtype}
\usepackage[strict,english=british]{csquotes}

\usepackage{amsmath,amssymb,amsthm}
\usepackage{mathtools}

\usepackage{algorithm}
\usepackage[noEnd=true,indLines=false]{algpseudocodex}

\usepackage[
	backend=biber,
	style=numeric-comp,
	sorting=nyt,
	sortcites=true,
	maxnames=5,
	giveninits=true,
	date=year,
	useprefix=true
]{biblatex}
\DeclareFieldFormat{eid}{no.~#1}
\usepackage{enumitem}

\usepackage[nobiblatex]{xurl}
\usepackage[breaklinks=true,pdflang=en-GB,colorlinks=true]{hyperref}
\colorlet{citecolor}{green!75!black}
\colorlet{linkcolor}{red!75!black}
\colorlet{urlcolor}{blue!75!black}
\hypersetup{citecolor=citecolor,linkcolor=linkcolor,urlcolor=urlcolor}
\usepackage{zref-clever}
\zcsetup{cap,nameinlink=false}
\usepackage{keytheorems}

\usepackage{orcidlink}
\usepackage{ellipsis}

\newkeytheorem{theorem}[parent=section,style=plain]
\newkeytheorem{lemma}[sibling=theorem,style=plain]
\newkeytheorem{proposition}[sibling=theorem,style=plain]

\newkeytheorem{definition}[sibling=theorem,style=definition]
\newkeytheorem{algothm}[sibling=theorem,style=definition,name=Algorithm,refname={algorithm,algorithms},Refname={Algorithm,Algorithms}]
\newkeytheorem{problem}[sibling=theorem,style=definition,refname={problem,problems},Refname={Problem,Problems}]

\newkeytheorem{remark}[sibling=theorem,style=remark]

\newcommand{\set}[2]{\ensuremath{\{\,{#1}\mid{#2}\,\}}}

\newcommand{\NN}{\mathbb{N}}

\DeclareMathOperator{\Hom}{Hom}
\DeclareMathOperator{\Aut}{Aut}
\DeclareMathOperator{\Aff}{Aff}
\DeclareMathOperator{\id}{id}
\DeclareMathOperator{\Coin}{Coin}
\DeclareMathOperator{\im}{im}
\newcommand{\bvarphi}{{\bar{\varphi}}}
\newcommand{\bpsi}{{\bar{\psi}}}
\newcommand{\tvarphi}{{\tilde{\varphi}}}
\newcommand{\tpsi}{{\tilde{\psi}}}
\newcommand{\restr}[2]{\ensuremath{{\left.\kern-\nulldelimiterspace{#1}\right|}_{#2}}}

\newcommand{\normalsub}{\trianglelefteq}
\DeclareMathOperator{\Fitt}{Fitt} 
\newcommand{\ab}[1]{#1^{\textup{ab}}}

\newcommand{\ind}[2]{\ensuremath{[#1\mathbin{:}#2]}}

\newcommand{\R}{{\mathop{\kern0pt \mathcal{R}}}}

\title[Algorithms for twisted conjugacy in polycyclic-by-finite groups]{Algorithms for twisted conjugacy classes of polycyclic-by-finite groups II}
\author[S. Tertooy]{Sam Tertooy\ \orcidlink{0000-0002-5750-9153}}
\date{\today}
\address{KU Leuven, Kulak Kortrijk Campus\\
	E.~Sabbelaan 53\\
	8500 Kortrijk\\
	Belgium}
\email{\href{mailto:sam.tertooy@kuleuven.be}{sam.tertooy@kuleuven.be}}
\urladdr{\url{https://stertooy.github.io}}


\begin{document}

\begin{abstract}
We construct an algorithm that, given a pair of homomorphisms between polycyclic-by-finite groups, determines whether their Reidemeister number is finite, and if so returns a set of representatives of the twisted conjugacy classes. Moreover, we show how this algorithm can be applied to compute double cosets and orbits of affine actions.
\end{abstract}

\maketitle

	\begin{center}
	This is an Accepted Manuscript of an article published by Elsevier in Journal of Algebra on 15 Sep 2026, available online: \href{https://doi.org/q5jb}{https://doi.org/q5jb}.
\end{center}

\section{Introduction}

Let \(G\) and \(H\) be groups and let \(\varphi,\psi \colon H \to G\) be group homomorphisms. Then \(H\) acts on \(G\) (from the right) by
\[	G \times H \to G \colon (g,h) \mapsto (h\varphi)^{-1}g(h\psi). \]
When \(g_1,g_2 \in G\) belong to the same orbit under this action, we say that they are \emph{\((\varphi,\psi)\)-twisted conjugate}.
These orbits are called the \emph{\((\varphi,\psi)\)-twisted conjugacy classes}, or the \emph{Reidemeister classes} of the pair \((\varphi,\psi)\). The orbit of \(g \in G\) is denoted by \([g]_{\varphi,\psi}\), and the set of all orbits by \(\R[\varphi,\psi]\). The \emph{Reidemeister number} \(R(\varphi,\psi)\) is the number of orbits and is either a positive integer or infinity.

The notion of twisted conjugacy finds its origin in topological coincidence theory; we refer to \cite{gonc05-a} for a survey on the subject. If \(f\) and \(g\) are continuous maps between topological spaces \(X\) and \(Y\), then the Reidemeister number \(R(f_*,g_*)\) of the induced group homomorphisms \(f_*, g_* \colon \pi_1(X) \to \pi_1(Y)\) holds information on the least number of coincidence points of (the homotopy classes of) \(f\) and \(g\).

One of the key components in extracting information from the Reidemeister number is the ability to determine its finiteness. This has led to the systematic search for groups with the \(R_\infty\)-property, i.e.\@ groups for which every Reidemeister number \(R(\varphi,\id)\) is infinite for every automorphism \(\varphi\). Recently, this property has been studied for soluble arithmetic groups \cite{ls23-a,ls24-a}, right-angled Artin groups \cite{ds21-a,lw24-a}, braid groups \cite{dgo21-a,dgo25-a}, linear groups \cite{ms22-a,ms24-a}, and generalised Baumslag-Solitar groups \cite{sw23-a,ssv23-a}, among other families.

The algorithmic study of twisted conjugacy has mostly been focused around solving the twisted conjugacy problem, i.e.\@ deciding whether or not two elements of a group are twisted conjugate, and variants thereof. In the past few years, solutions were found for e.g.\@ Artin groups \cite{crow25-a,crow26-a,bv24-a}, direct products of free groups \cite{cd26-a,carv26-a}, and soluble Baumslag-Solitar groups \cite{mrv25-a,roy26-a}. Closely related to the present paper are the results obtained by Roman'kov for finitely generated metabelian and polycyclic groups \cite{roma10-a,roma16-a,roma21-a,rv09-a}.

In the current paper, we instead focus on the algorithmic computation of Reidemeister numbers and representatives of Reidemeister classes. We consider the following two problems:

\begin{problem}[manual-num=A]
    \label{prob:introproblem1}
    Given groups \(G, H\) and homomorphisms \(\varphi,\psi \colon H \to G\), compute the Reidemeister number \(R(\varphi,\psi)\).
\end{problem}

\begin{problem}[manual-num=B]
    \label{prob:introproblem2}
    Given groups \(G, H\) and homomorphisms \(\varphi,\psi \colon H \to G\) with Reidemeister number \(r < \infty\), find \(\{g_1, g_2, \ldots, g_r\} \subseteq G \) such that 
    \( G = [g_1]_{\varphi,\psi} \sqcup [g_2]_{\varphi,\psi} \sqcup \dots \sqcup [g_r]_{\varphi,\psi}\).
\end{problem}

Our primary goal is to construct an algorithm that solves these problems when both \(G\) and \(H\) are polycyclic-by-finite. In \cite{dt21-a}, we addressed the case \(G = H\). However, that solution does not extend straightforwardly to the case \(G \neq H\). Instead, we construct a new algorithm that solves a problem which generalises both \zcref{prob:introproblem1,prob:introproblem2}.

Consider the situation from \zcref{prob:introproblem1} and let \(N\) be a normal subgroup of \(G\). We denote the set of \((\varphi,\psi)\)-twisted conjugacy classes that intersect \(N\) non-trivially by \({\R}_N[\varphi,\psi]\), and the number of such classes by \(R_N(\varphi,\psi)\). We construct the following algorithm.

\begin{algothm}[manual-num=A]
    \label{alg:mainalgA}
    There exists an algorithm that, given two polycyclic-by-finite groups \(G\) and \(H\), two homomorphisms \(\varphi,\psi \colon H \to G\), and a normal subgroup \(N \normalsub G\), computes \(R_N(\varphi,\psi)\), and if it is finite, finds \(\{n_1, n_2, \ldots, n_r\} \subseteq N\) such that \(N \subseteq [n_1]_{\varphi,\psi} \sqcup [n_2]_{\varphi,\psi} \sqcup \dots \sqcup [n_r]_{\varphi,\psi}\).
\end{algothm}

By taking \(N = G\), it is clear that this algorithm provides a solution to \zcref{prob:introproblem1,prob:introproblem2}. An implementation (for polycyclic groups only) is available in the \textsf{GAP} \cite{gap25-a} package \texttt{TwistedConjugacy} \cite{tert26-a}.

This paper is structured as follows. \zcref[cap]{sec:prelims} introduces the necessary preliminaries on polycyclic-by-finite groups and twisted conjugacy. In \zcref{sec:setup,sec:nilpotent,sec:metabelian,sec:abcd,sec:general} we construct \zcref{alg:mainalgA} by gradually considering more general cases.
Finally, \zcref{sec:doublecosets,sec:affineactions} illustrate how algorithms for twisted conjugacy can be applied to computations with double cosets and affine actions, respectively.

\section{Preliminaries}\label{sec:prelims}
First, let us fix some notation. We work exclusively with right group actions, and similarly always apply maps from the right. If \(f\colon X \to Y\) is a map and \(x \in X\), then the image of \(x\) under \(f\) is denoted by \(xf\). The composition \(fg\) of maps \(f\) and \(g\) is thus defined by \(x(fg) = (xf)g\). If \(x,y\) are elements of a group, then we use exponents to denote the conjugation action (\(x^y \coloneq y^{-1}xy\)) and square brackets to denote the commutator (\([x,y] \coloneq x^{-1}y^{-1}xy\)). If \(G\) is a group and \(g \in G\), then \(\iota_g\) denotes the inner automorphism of \(G\) given by \(x \mapsto x^g\). The semi-direct product of a group \(H\) acting on a group \(G\) will be denoted by \(H \ltimes G\), and its multiplication is defined as \((h_1,g_1)(h_2,g_2) = (h_1h_2,g_1^{h_2}g_2)\).

In the remainder of this section, we briefly cover some of the properties (algorithmic or otherwise) of polycyclic-by-finite groups that we use in the sequel. For a treatise on polycyclic (and by extension polycyclic-by-finite) groups, we refer to \cite{sega83-a}.

\begin{proposition}
    \label{prop:pbfproperties}
    Let \(G\) be a polycyclic-by-finite group. Then
    \begin{enumerate}[label=(\roman*)]
        \item \(G\) is nilpotent-by-abelian-by-finite,
        \item every subgroup of \(G\) is finitely generated,
        \item every ascending subgroup series in \(G\) stabilises.
    \end{enumerate}
\end{proposition}

In particular, if \(Z_i(G)\) denotes the \(i\)-th term of the upper central series of a polycyclic-by-finite group \(G\), i.e.\@
\[Z_0(G) = 1, \quad \frac{Z_{i+1}(G)}{Z_i(G)} = Z\left(\frac{G}{Z_i(G)}\right),  \]
then there exists a \(k \in \NN\) such that \(Z_k(G) = Z_{k+1}(G)\).

We will also need the following result on the separability of twisted conjugacy classes in polycyclic-by-finite groups, which was obtained in \cite[Thm.~5.1]{tert25-b}.

\begin{theorem}
    \label{thm:sep}
    Let \(G\) be a polycyclic-by-finite group and \(H\) any group. Let \(\varphi,\psi \in \Hom(H,G)\) and \(g_1,g_2 \in G\). If \(g_1 \notin [g_2]_{\varphi,\psi}\), then there exists a finite index normal subgroup \(N\) of \(G\) such that \(g_1 \notin [g_2]_{\varphi,\psi} N\), or, equivalently, such that \([g_1p]_{\varphi p,\psi p} \neq [g_2p]_{\varphi p,\psi p}\) with \(p \colon G \to G/N\) the natural projection.
\end{theorem}

We refer to \cite{bcr91-a,eick01-a} for a comprehensive treatment of the algorithmic theory of polycyclic(-by-finite) groups. In what follows, when we say ``let \(G\) be a polycyclic-by-finite group'', we will always mean that \(G\) is given by a finite presentation. For practical purposes, if \(G\) is polycyclic, the calculation of a so-called \emph{consistent polycyclic presentation} is paramount. If \(G\) is not polycyclic, then one can calculate a consistent polycyclic presentation for a polycyclic normal subgroup, and use the techniques presented in \cite{sh17-a} for efficient computation. 

To summarise some of the results we require, let \(G,H\) be polycyclic-by-finite groups, \(K,L\) subgroups and \(N\) a normal subgroup of \(G\), and \(\varphi \in \Hom(H,G)\). Then there exist algorithms that compute or construct the following:
\begin{itemize}
    \item the Fitting subgroup \(\Fitt(G)\) (i.e.\@ the maximal nilpotent normal subgroup),
    \item the derived subgroup \(G'\),
    \item the centre \(Z(G)\),
    \item a nilpotent-by-abelian finite index normal subgroup of \(G\),
    \item the index \(\ind{G}{K}\),
    \item the intersection \(K \cap L\),
    \item the product \(NK\),
    \item the image \(H\varphi\),
    \item the preimage \(K\varphi^{-1}\).
\end{itemize}
Moreover, there is an algorithm to recursively enumerate the finite index normal subgroups of a polycyclic-by-finite group \cite[Cor.~2.7]{bcr91-a}.

At the heart of the algorithmic study of twisted conjugacy are the \emph{twisted conjugacy (decision) problem} 
 and the \emph{twisted conjugacy search problem}.
\begin{problem}
    \label{prob:tcp}
    Let \(G, H\) be groups, let \(\varphi,\psi \in \Hom(H,G)\), and let \(g_1,g_2 \in G\). Decide whether \(g_1\) and \(g_2\) are \((\varphi,\psi)\)-twisted conjugate.
\end{problem}

\begin{problem}
    \label{prob:tcsp}
    Let \(G, H\) be groups, let \(\varphi,\psi \in \Hom(H,G)\), and let \(g_1,g_2 \in G\) be \((\varphi,\psi)\)-twisted conjugate. Find \(h \in H\) such that \((h\varphi)^{-1}g_1(h\psi) = g_2\).
\end{problem}

The algorithm below provides a solution to both problems.
\begin{algothm}
    \label{alg:tcp}
    There exists an algorithm that, given two polycyclic-by-finite groups \(G\) and \(H\), two homomorphisms \(\varphi,\psi \colon H \to G\), and \(g_1,g_2 \in G\), determines whether \(g_1\) and \(g_2\) are \((\varphi,\psi)\)-twisted conjugate, and if they are, finds \(h \in H\) such that \((h\varphi)^{-1}g_1(h\psi) = g_2\).
\end{algothm}
\begin{proof}
    The so-called \emph{local-global method}, also known as McKinsey's algorithm \cite{mcki43-a}, can be applied here. We start two procedures in parallel.
    
    The first procedure iterates over all elements \(h\) of \(H\). For each \(h\), it calculates \((h\varphi)^{-1}g_1(h\psi)\) and tests whether this element is equal to \(g_2\). If it is, then we halt both procedures and the algorithm terminates.
    
    The second procedure iterates over all finite index normal subgroups \(M\) of \(G\). For each \(M\), it calculates subsequently \(N \coloneq M\varphi^{-1} \cap M\psi^{-1}\), the finite quotients \(\bar{G} \coloneq G/M\) and \(\bar{H} \coloneq H/N\), the projections \(p \colon G \to \bar{G}\) and \(q\colon H \to \bar{H}\), and the induced homomorphisms \(\bvarphi\) and \(\bpsi\) defined by \(\varphi p = q \bvarphi\) and \(\psi p = q \bpsi\), respectively.
    Finally, it tests whether \(g_1p\) and \(g_2p\) are \((\bvarphi,\bpsi)\)-twisted conjugate. If they are not, then \(g_1\) and \(g_2\) are not \((\varphi,\psi)\)-twisted conjugate either, so we halt both procedures and the algorithm terminates. 
    
    If \(g_1\) and \(g_2\) are \((\varphi,\psi)\)-twisted conjugate, the first procedure is guaranteed to eventually find a suitable \(h\). If they are not twisted conjugate, then by \zcref{thm:sep} the second procedure is guaranteed to eventually find a suitable \(M\). Therefore this algorithm will always terminate.
\end{proof}

While straightforward to construct, the algorithm proposed in the above proof is certainly not expected to be efficient. A more practical algorithm is described in \cite[Sec.~7]{roma21-a} and \cite[Sec.~5.4.6]{bkl20-a}, and has been implemented (for polycyclic groups) in the \texttt{TwistedConjugacy} package.

Related to the twisted conjugacy problem is the \emph{coincidence problem}, also known as the \emph{equaliser problem}. This problem concerns the \emph{coincidence group}
\[ \Coin(\varphi,\psi) \coloneq  \set{ h \in H }{h\varphi = h\psi},\]
where \(\varphi,\psi \in \Hom(H,G)\) for certain groups \(G,H\).

\begin{problem}
    \label{prob:coin}
    Let \(G\) and \(H\) be groups and let \(\varphi,\psi \in \Hom(H,G)\). Find a presentation for \(\Coin(\varphi,\psi)\).
\end{problem}

In the case of polycyclic-by-finite groups, it is sufficient to find a generating set of \(\Coin(\varphi,\psi)\). A complete presentation of \(\Coin(\varphi,\psi)\) can then be calculated from the presentation of \(H\), if desired, e.g.\@ using the algorithm described in \cite[Thm.~3.4]{bcr91-a}.

\begin{algothm}
    \label{alg:coin}
    There exists an algorithm that finds a generating set for the coincidence group of a pair of homomorphisms between polycyclic-by-finite groups.
\end{algothm}
\begin{proof}
    Consider the group homomorphism
    \[\varphi \times \psi \colon H \to G \times G \colon h \mapsto (h\varphi,h\psi),\]
    and the preimage \(\Delta_G(\varphi \times \psi)^{-1}\) of the diagonal \(\Delta_G = \{\,(g,g) \in G \times G \mid g \in G \,\}\). It is easy to see that the coincidence group \(\Coin(\varphi,\psi)\) coincides with this preimage.
\end{proof}

Again, more practical algorithms exist; these can be found in \cite[Sec.~7]{roma21-a} and \cite[Sec.~5.4.6.2]{bkl20-a} and have been implemented in the \texttt{TwistedConjugacy} package.

\begin{remark}
    \label{rem:tcpcoingeneralisation}
Neither \zcref{alg:tcp} nor \zcref{alg:coin} directly exploits the structure of polycyclic-by-finite groups. Rather, we used that twisted conjugacy classes are separable (in the first algorithm), that the class of polycyclic-by-finite groups is closed under taking direct products (in the second algorithm), and that we are able to compute certain subgroups, intersections, (pre)images, etc.\@ as discussed below \zcref{thm:sep}.
Thus, we could apply these algorithms to a group which fulfils the requirements on separability and/or direct products, and for which methods to compute the subgroups, (pre)images, etc.\@ are known to exist.
\end{remark}

\section{Setup and tools}
\label{sec:setup}
In this section, we introduce the framework and tools we use to construct \zcref{alg:mainalgA}. To start off, we introduce a slight generalisation of the coincidence group. Consider the situation from \zcref{alg:coin} and let \(N\) be a normal subgroup of \(G\). We define
\[ \Coin_N(\varphi,\psi) \coloneq \set{ h \in H }{(h\varphi)^{-1}(h\psi) \in N}.\]
Let \(p\colon G \to G/N\) be the natural projection and set \(\bvarphi \coloneq \varphi p\), \(\bpsi \coloneq \psi p\). Then \(\Coin_N(\varphi,\psi) = \Coin(\bvarphi,\bpsi)\). Thus, every such group can also be calculated using \zcref{alg:coin}.

The first step in \zcref{alg:mainalgA} is replacing \(H\) by \(\Coin_N(\varphi,\psi)\), and subsequently restricting \(\varphi\) and \(\psi\) to this subgroup. The resulting action has the same orbits, but can be restricted to an action of \(\Coin_N(\varphi,\psi)\) on \(N\).
\begin{definition}
    \label{def:stdquintuple}
    A quintuple \((G,H,\varphi,\psi,N)\) where
    \begin{itemize}
        \item \(G\) and \(H\) are polycyclic-by-finite groups,
        \item \(\varphi\) and \(\psi\) are group homomorphisms \(H \to G\),
        \item \(N\) is a normal subgroup of \(G\) such that \(H = \Coin_N(\varphi,\psi)\),
    \end{itemize}
    will be called a \emph{standard quintuple}.
\end{definition}

The construction of \zcref{alg:mainalgA} now reduces to an algorithm which solves the \zcref*[noref,nocap]{prob:mainproblem} below.

\begin{problem}[manual-num=C]
    \label{prob:mainproblem}
    Given a standard quintuple \((G,H,\varphi,\psi,N)\), compute \(R_N(\varphi,\psi)\), and if it is finite, find \(\{n_1, n_2, \ldots, n_r\} \subseteq N \) such that \(N = [n_1]_{\varphi,\psi} \sqcup [n_2]_{\varphi,\psi} \sqcup \dots \sqcup [n_r]_{\varphi,\psi}\).
\end{problem}

We provide two ``auxiliary'' algorithms that will prove to be indispensable in the coming sections. The first algorithm reduces the calculation of twisted conjugacy classes in \(N \normalsub G\) to that of twisted conjugacy classes in \(N \cap K \normalsub G\) and in \(NK/K \normalsub G/K\), for a normal subgroup \(K \normalsub G\). The underlying idea is described in the next \zcref*[noref,nocap]{lem:reidclassdecomp}.

\begin{lemma}
    \label{lem:reidclassdecomp}
    Let \((G,H,\varphi,\psi,N)\) be a standard quintuple and let \(K \normalsub G\). Let \(p \colon G \to G/K\) be the natural projection, set \(M \coloneq N \cap K\) and define
    \begin{align*}
        &\bar{G} \coloneq Gp, &&\bvarphi \coloneq \varphi p,\\
        &\bar{N} \coloneq Np, &&\bpsi \coloneq \psi p.
    \end{align*}
    Then, for every \(n \in N\), define
    \begin{gather*}
        C_n \coloneq \Coin_K(\varphi\iota_{n},\psi),\\
        \varphi_n \colon C_n \to G \colon c \mapsto c(\varphi\iota_{n}),
        \qquad \psi_n \colon C_n \to G \colon c \mapsto c\psi.
    \end{gather*}
    Let \(n_1, n_2, \dotsc \in N\) such that \(\bar{N} = \bigsqcup_{i} [\bar{n}_i]_{\bvarphi,\bpsi}\), where \(\bar{n}_i \coloneq n_ip\), and for each \(n_i\), let
    \(m_{n_i1}, m_{n_i2}, \dotsc \in M\) such that \( M = \bigsqcup_j [m_{n_ij}]_{\varphi_{n_i},\psi_{n_i}}.\)
    Then \(N\) is the (disjoint) union of the following \((\varphi,\psi)\)-twisted conjugacy classes:
    \[ N = \bigsqcup_{i}\bigsqcup_{j} [n_im_{n_ij}]_{\varphi,\psi}.\]
    In particular, \(R_N(\varphi,\psi)\) is finite if and only if \(R_{\bar{N}}(\bvarphi,\bpsi)\) is finite and \(R_{M}(\varphi_{n_i},\psi_{n_i})\) is finite for every \(n_i\).
\end{lemma}
\begin{proof}
    It suffices to prove that every \(n \in N\) belongs to exactly one twisted conjugacy class \([n_im_{n_ij}]_{\varphi,\psi}\). Let \(n \in N\) and consider \(\bar{n} \coloneq np\). 
    Let \(i \in \NN\) be such that \([\bar{n}]_{\bvarphi,\bpsi} = [\bar{n}_i]_{\bvarphi,\bpsi}\). There exists an \(h_1 \in H\) such that
    \[\bar{n} = (h_1\bvarphi)^{-1}\bar{n}_i(h_1\bpsi).\]
    Lifting this back to \(G\), there exists an \(m \in M\) for which
    \[
        n =  (h_1\varphi)^{-1}n_im(h_1\psi) =  n_i(h_1\varphi\iota_{n_i})^{-1}m(h_1\psi).
    \]
    Now let \(j \in \NN\) such that \([m]_{\varphi_i,\psi_i} = [m_{n_ij}]_{\varphi_i,\psi_i}\). There exists an \(h_2 \in H\) such that
    \[
        m = (h_2\varphi_i)^{-1}m_{n_ij}(h_2\psi_i) = (h_2\varphi\iota_{n_i})^{-1} m_{n_ij}(h_2\psi).
    \]
    Combining the previous two equations and setting \(h \coloneq h_2h_1\), we obtain that
    \begin{equation*}
        n
        = n_i  (h\varphi\iota_{n_i})^{-1}m_{n_ij}(h\psi)
        = (h\varphi)^{-1}n_im_{n_ij} (h\psi),
    \end{equation*}
    hence \(n \in [n_im_{n_ij}]_{\varphi,\psi}\) for some \(i, j \in \NN\). We now confirm that the union is disjoint. Consider representatives \(n_i m_{n_ij}\) and \(n_k m_{n_kl}\). Then
    \begin{align*}
    [n_i m_{n_ij}]_{\varphi,\psi} = [n_k m_{n_kl}]_{\varphi,\psi}
    \iff& \exists h \in H \colon n_i m_{n_ij} = (h\varphi)^{-1}n_k m_{n_kl}(h\psi)\\
    \implies& \exists h \in H \colon \bar{n}_i = (h\bvarphi)^{-1}\bar{n}_k(h\bpsi)\\
    \iff& [\bar{n}_i]_{\bvarphi,\bpsi} = [\bar{n}_k]_{\bvarphi,\bpsi}\\
    \iff& n_i = n_k.
    \end{align*}
    By contraposition, if \(i \neq k\), then the classes \([n_i m_{n_ij}]_{\varphi,\psi}\) and \([n_k m_{n_kl}]_{\varphi,\psi}\) are distinct.
    If \(i = k\), then
    \begin{align*}
    [n_i m_{n_ij}]_{\varphi,\psi} = [n_i m_{n_il}]_{\varphi,\psi}
    \iff& \exists h \in H \colon n_i m_{n_ij} = (h\varphi)^{-1}n_i m_{n_il}(h\psi)\\
    \iff& \exists h \in H \colon m_{n_ij} = (h\varphi\iota_{n_i})^{-1}m_{n_il}(h\psi)\\
    \iff& \exists h \in C_{n_i} \colon m_{n_ij} = (h\varphi_{n_i})^{-1}m_{n_il}(h\psi_{n_i})\\
    \iff& [m_{n_ij}]_{\varphi_{n_i},\psi_{n_i}} = [m_{n_il}]_{\varphi_{n_i},\psi_{n_i}}\\
    \iff& m_{n_ij} = m_{n_il}.
    \end{align*}
    We conclude that \(\bigsqcup_{i}\bigsqcup_{j} [n_im_{n_ij}]_{\varphi,\psi}\) is a disjoint union.
\end{proof}

\begin{algothm}
    \label{alg:decomp}
    Let \((G,H,\varphi,\psi,N)\) be a standard quintuple and let \(K \normalsub G\). Using the definitions from \zcref{lem:reidclassdecomp}, suppose that there exists an algorithm that solves \zcref{prob:mainproblem} for  \((\bar{G},H,\bvarphi,\bpsi,\bar{N})\), and for every \(n \in N\) there exists an algorithm that does the same for \((G,C_n,\varphi_n,\psi_n,M)\). Then there exists an algorithm that solves \zcref{prob:mainproblem} for \((G,H,\varphi,\psi,N)\).
\end{algothm}

One limitation of this \zcref*[noref,nocap]{alg:decomp} is that the group \(G\) remains unchanged when we pass from \(N\) to \(M = N \cap K\). By contrast, the second \zcref*[noref,nocap]{alg:byfinite} gives us a way to pass from \(G\) to a finite index normal subgroup of \(G\). 

\begin{algothm}
    \label{alg:byfinite}
    Let \((G,H,\varphi,\psi,N)\) be a standard quintuple and let \(K \normalsub G\) with \(\ind{G}{K} < \infty\). Suppose that there exists an algorithm that solves \zcref{prob:mainproblem} for every standard quintuple of the form \((K,L,\lambda,\mu,M)\).
    There exists an algorithm that solves \zcref{prob:mainproblem} for \((G,H,\varphi,\psi,N)\).
\end{algothm}
\begin{proof}
    We consider two cases. First, suppose that \(N\) is contained in \(K\). Define
    \begin{gather*}
        L \coloneq  K\varphi^{-1} \cap K\psi^{-1},\\
\lambda \colon L \to K \colon l \mapsto l\varphi,\qquad
\mu \colon L \to K \colon l \mapsto l\psi,
    \end{gather*}
    and consider the (surjective) map
    \[ \pi \colon \R_N[\lambda,\mu] \to \R_N[\varphi,\psi] \colon [n]_{\lambda,\mu} \mapsto [n]_{\varphi,\psi}.\]
    Since \(K\) is a finite index normal subgroup of \(G\), we deduce that \(L\) is a finite index normal subgroup of \(H\), and thus we can pick a transversal \(\{h_1, \ldots, h_r\}\) of \(L\) in \(H\). Now suppose that \(n,n' \in N\) are such that \([n]_{\varphi,\psi} = [n']_{\varphi,\psi}\), i.e.\@ there exists some \(h \in H\) such that \(n' = (h\varphi)^{-1}n(h\psi)\). For some (unique) \(h_i\) and some \(l \in L\), we have \(h = h_il\). Now
    \[n' = (h\varphi)^{-1}n(h\psi) = (l\varphi)^{-1}(h_i\varphi)^{-1}n(h_i\psi)(l\psi)= (l\lambda)^{-1}(h_i\varphi)^{-1}n(h_i\psi)(l\mu),\]
    which implies that
    \[[n']_{\lambda,\mu} = [(h_i\varphi)^{-1}n(h_i\psi)]_{\lambda,\mu}.\]
    Hence for every \(n \in N\), the preimage \([n]_{\varphi,\psi}\pi^{-1}\) is finite, and in particular \(R_N(\lambda,\mu) < \infty\) if and only if \(R_N(\varphi,\psi) < \infty\). If \(R_N(\lambda,\mu)\) is infinite, then so is \(R_N(\varphi,\psi)\) and the algorithm finishes. Otherwise, we obtain \(\{n_1, n_2, \ldots, n_t\} \subseteq N\) such that
    \(N = [n_1]_{\lambda,\mu} \sqcup [n_2]_{\lambda,\mu} \sqcup \cdots \sqcup [n_t]_{\lambda,\mu}\).
    
    Applying \zcref{alg:tcp} to pairs of \(n_i\)'s, we can reduce this to a finite set \(\{n_{i_1}, n_{i_2}, \ldots, n_{i_s}\} \) such that
    \(N = [n_{i_1}]_{\varphi,\psi} \sqcup [n_{i_2}]_{\varphi,\psi} \sqcup \cdots \sqcup [n_{i_s}]_{\varphi,\psi}\),
    which finishes the first case.
    
    We now move on to the second case: \(N\) is not contained in \(K\). We apply \zcref{alg:decomp} for the normal subgroup \(K\). Indeed, \(\bar{G}\) is finite, which poses no problem, and for each of the quintuples \((G,C_{n_i},\varphi_{n_i},\psi_{n_i},M)\) we note that \(M \leq K\), hence we can repeat the steps described in the first case.
\end{proof}

\section{Nilpotent-by-finite groups}
\label{sec:nilpotent}
The main result of this section is the existence of the \zcref*[noref,nocap]{alg:nilpotentbyfinite} below.
\begin{algothm}
    \label{alg:nilpotentbyfinite}
    There exists an algorithm that solves \zcref{prob:mainproblem} for every standard quintuple \((G,H,\varphi,\psi,N)\) where \(G\) is nilpotent-by-finite.
\end{algothm}

Rather than providing one large algorithm immediately, we split this up into several ``sub-algorithms''. A first step is to consider the case where \(N\) is a central subgroup of \(G\).

\begin{algothm}
    \label{alg:central}
    There exists an algorithm that solves \zcref{prob:mainproblem} for every standard quintuple \((G,H,\varphi,\psi,N)\) where \(N\) is central.
\end{algothm}
\begin{proof}
    We start by constructing the group homomorphism given by
    \[\delta\colon H \to N \colon h \mapsto (h\varphi)^{-1}(h\psi).\]
    Let \(p\colon N \to N/H\delta\) be the natural projection and set \(\bar{N} \coloneq Np\). Since \(N\) is central, for every \(n_1, n_2 \in N\),
    \begin{align*}
    [n_1]_{\varphi,\psi} = [n_2]_{\varphi,\psi}
    &\iff \exists h \in H \colon n_1 = (h\varphi)^{-1}n_2(h\psi)\\
    &\iff \exists h \in H \colon n_1 = n_2(h\varphi)^{-1}(h\psi)\\
    &\iff \exists h \in H \colon n_1 = n_2(h\delta)\\
    &\iff n_1p = n_2p.
    \end{align*}
    It follows that the map
    \[\pi \colon {\R}_N[\varphi,\psi] \to \bar{N} \colon [n]_{\varphi,\psi} \mapsto np\]
    is a bijection. 
    Thus, if \(\bar{N}\) is infinite, so is \(R_N(\varphi,\psi)\). Otherwise, denote the elements of \(\bar{N}\) by \(\bar{n}_1, \ldots, \bar{n}_k\). For each \(i \in \{1,\ldots,k\}\), pick \(n_i \in \bar{n}_ip^{-1}\). Then \(\{n_1, \ldots, n_k\}\) satisfies \(N = \bigsqcup_{i=1}^k [n_i]_{\varphi,\psi},\)
    so we have solved \zcref{prob:mainproblem}.
\end{proof}

If \(G\) is a nilpotent group, then \(Z_k(G) = G\) for some \(k \in \NN\). This hints at the possibility of tackling the case where \(G\) is nilpotent by inductively applying the previous algorithm.

\begin{algothm}
    \label{alg:nilpotent}
    There exists an algorithm that solves \zcref{prob:mainproblem} for every standard quintuple \((G,H,\varphi,\psi,N)\) where \(G\) is nilpotent.
\end{algothm}

\begin{proof}
    Let \((G,H,\varphi,\psi,N)\) be a standard quintuple with \(G\) nilpotent of class \(c\). We prove this by induction on \(c\). If \(c = 0\), then \(G\) is trivial and therefore \(N = [1]_{\varphi,\psi}\).
    
    Now suppose that \(c > 0\). We can apply \zcref{alg:decomp} for the normal subgroup \(K \coloneq Z(G)\). Indeed, if we consider the quintuple \((\bar{G},H,\bvarphi,\bpsi,\bar{N})\), then \(\bar{G}\) is nilpotent of class \(c-1\), so the required algorithm exists by induction. For every \(n \in N\) the quintuple \((G,C_n,\varphi_n,\psi_n,M)\) has \(M \leq Z(G)\), so we can apply \zcref{alg:central}.
\end{proof}

Finally, the general case of \(G\) being nilpotent-by-finite now follows easily by combining previously obtained algorithms.

\begin{proof}[Proof of \zcref{alg:nilpotentbyfinite}]
    Let \(K\) be the Fitting subgroup of \(G\). Since \zcref{alg:nilpotent} can solve \zcref{prob:mainproblem} for every standard quintuple of the form \((K,L,\lambda,\mu,M)\), the result follows from \zcref{alg:byfinite}.
\end{proof}

\section{Metabelian groups}
\label{sec:metabelian}

As in the previous section, the main result here is an algorithm solving \zcref{prob:mainproblem} when some additional conditions are placed on the standard quintuple. In particular, \(G\) will be metabelian throughout this section.

\begin{algothm}
    \label{alg:metabelian}
    There exists an algorithm that solves \zcref{prob:mainproblem} for every standard quintuple \((G,H,\varphi,\psi,A)\) where
    \begin{enumerate}
        \item \(H\) is abelian,
        \item \(A\) is abelian,
        \item \(G = AH\varphi\).
    \end{enumerate}
    Such a quintuple will be called a \emph{metabelian quintuple}. 
\end{algothm}

We construct this algorithm with the help of some theoretical results. This will require a very brief and basic introduction to group derivations; broader expositions can be found in most standard works on
group cohomology (e.g.\@ \cite[Sec.~IV.2]{brow82-a}). We denote the (right) action of a group \(Q\) on a group \(A\) by exponents, i.e.\@ the action of \(Q\) on \(A\) is given by the map
\[A \times Q \to A \colon (a,q) \mapsto a^q\]
Let us also introduce the notation \([a,q] \coloneq a^{-1}a^q\). We remark that both notations introduced here coincide with those of the previous section if \(A\) and \(Q\) are subgroups of a common supergroup and \(Q\) acts on \(A\) by conjugation. To remain consistent with the remainder of this paper, we use multiplicative notation even when the groups being discussed are abelian.

\begin{definition}
    \label{def:module}
    Let \(Q\) be a group acting on an abelian group \(A\) such that for all \(a_1,a_2 \in A\) and all \(q \in Q\) we have
    \((a_1a_2)^q = a_1^q a_2^q\). Then \(A\) is called a \emph{\(Q\)-module}.
\end{definition}

An equivalent way to define a module is to say that \(Q\) acts \emph{via automorphisms} on \(A\), i.e.\@ there is a homomorphism \(Q \to \Aut(A) \colon q \mapsto \lambda_q\) such that \(a^q = a\lambda_q\).

\begin{definition}
    \label{def:derivation}
    Let \(Q\) be a group and \(A\) a \(Q\)-module. A map \( \delta \colon Q \to A \) is called a \emph{derivation} (or \emph{crossed homomorphism}) if \((q_1q_2)\delta = (q_1\delta)^{q_2}(q_2\delta)\) for all \(q_1,q_2 \in Q\).
\end{definition}

When it may not be clear from the context what the \(Q\)-module structure on \(A\) is, we will say that a map \(\delta \colon Q \to A\) is a derivation \emph{with respect to} a particular action of \(Q\) on \(A\) or a homomorphism \(Q \to \Aut(A)\).

\begin{definition}
    \label{def:cohomology}
    Let \(Q\) be a group and let \(A\) be a \(Q\)-module. By \(A^Q\) we denote the set of all \(Q\)-invariant elements of \(A\), i.e.\@
    \[ A^Q \coloneq \set{a \in A}{\forall q \in Q: [a,q] = 1}.\]
\end{definition}

We now have the necessary background to state the following \zcref[noref,nocap]{thm:surjectivederivation}, which can be extracted from part (iii) in the proof of \cite[Thm.~B]{robi02-a}.
\begin{theorem}
    \label{thm:surjectivederivation}
    Let \(A\) and \(Q\) be finitely generated abelian groups such that \(A\) is a \(Q\)-module with \(A^Q = 1\).
    If a surjective derivation \(\delta \colon Q \to A\) exists, then \(A\) is finite.
\end{theorem}

Having finished our introduction to group derivations and with the above \zcref*[noref,nocap]{thm:surjectivederivation} obtained, we are now ready to tackle metabelian quintuples. We prove three \zcref*[noref,nocap]{lem:ANquint1,lem:ANquint2,lem:ANquint3}, each building on the previous.

\begin{lemma}
    \label{lem:ANquint1}
    Let \((G,H,\varphi,\psi,A)\) be a metabelian quintuple such that \(R_A(\varphi,\psi) = 1\) and \(Z(G) = 1\). Then \(A\) is finite.
\end{lemma}
\begin{proof}
    The group \(H\) acts on \(A\) via 
    \[A \times H \to A \colon (a,h)\mapsto a^h \coloneq (h\varphi)^{-1}a(h\varphi),\]
    so \(A\) is an \(H\)-module. The map \(\delta\) defined by
    \[\delta \colon H \to A \colon h \mapsto (h\varphi)^{-1}(h\psi)\]
    is a derivation, and it is surjective since \(A = [1]_{\varphi,\psi} = H\delta\). Now suppose \(a \in A^H\), i.e.\@ \(1 = [a,h] = [a,h\varphi]\) for all \(h \in H\). Since \(G = AH\varphi\), it then follows that \(a \in Z(G)\) and therefore \(a = 1\), so \(A^H\) is trivial. By \zcref{thm:surjectivederivation}, \(A\) is finite.
\end{proof}

\begin{lemma}
    \label{lem:ANquint2}
    Let \((G,H,\varphi,\psi,A)\) be a metabelian quintuple such that \(R_A(\varphi,\psi) = 1\). Then \(G\) is nilpotent-by-finite.
\end{lemma}
\begin{proof}
    As \(G\) is polycyclic, its upper central series eventually stabilises, so \(Z_k(G) = Z_{k+1}(G)\) for some \(k \in \NN\). Let \(p \colon G \to G/Z_k(G)\) be the natural projection and set \(\bar{G} \coloneq Gp\), \(\bar{A} \coloneq Ap\), \(\bvarphi \coloneq \varphi p\) and \(\bpsi \coloneq \psi p\).
    
    Then \((\bar{G},H,\bvarphi,\bpsi,\bar{A})\) is a metabelian quintuple that satisfies the conditions of \zcref{lem:ANquint1}, hence \(\bar{A}\) is finite and \(\bar{G}\) is therefore finite-by-abelian. It follows from \cite[Lem.~6.3]{tert25-b} that \(\bar{G}\) is abelian-by-finite. Now let \(\bar{B}\) be an abelian finite index normal subgroup of \(\bar{G}\), set \(B \coloneq \bar{B}p^{-1}\) and set \(\bar{Z} \coloneq Z_k(B)p\). Since \(Z_k(G) \leq B\), we deduce that \(Z_i(G) \leq Z_i(B)\) for all \(i \in \NN\). The Third Isomorphism Theorem implies that
    \[ \frac{B}{Z_k(B)} \cong \frac{B/Z_k(G)}{Z_k(B) / Z_k(G)} = \frac{\bar{B}}{\bar{Z}}.\]
    Since \(B/Z_k(B)\) is isomorphic to a quotient of the abelian group \(\bar{B}\), it is itself abelian. But then \(Z_{k+1}(B) = B\), hence \(B\) is nilpotent. As \(B\) has finite index in \(G\), the latter is nilpotent-by-finite.
\end{proof}

\begin{lemma}
    \label{lem:ANquint3}
    Let \((G,H,\varphi,\psi,A)\) be a metabelian quintuple such that \(R_A(\varphi,\psi) < \infty\). Then \(G\) is nilpotent-by-finite.
\end{lemma}
\begin{proof}
    Let \(a_1 = 1, a_2, \ldots, a_n \in A\) be such that \( A = \bigsqcup_{i=1}^n [a_i]_{\varphi,\psi}\).
    Invoking \zcref{thm:sep}, for each \(i \in \{2,\dotsc,n\}\), there exists a finite index normal subgroup \(N_i\) of \(G\) such that \(a_i \notin [1]_{\varphi,\psi}N_i\). Define \(N \coloneq \bigcap_{i=2}^n N_i\), now \(B \coloneq A \cap N\) has finite index in \(A\), it is normal in \(G\) and is a subset of \([1]_{\varphi,\psi}\). Next, we set \(C \coloneq \Coin_N(\varphi,\psi)\), \(K \coloneq B C\varphi\), and we define \(\lambda,\mu \in \Hom(C,K)\) by
    \[
        \lambda \colon C \to K \colon c \mapsto c\varphi,\quad
        \mu \colon C \to K \colon c \mapsto c\psi.
    \]
    It is straightforward to verify that the quintuple \((K,C,\lambda,\mu,B)\) satisfies the conditions of \zcref{lem:ANquint2}, so \(K\) is nilpotent-by-finite. Let \(p \colon G \to G/N\) be the natural projection. Then
    \[ \ker(\varphi p) \cap \ker(\psi p) \leq \Coin(\varphi p, \psi p) = C.\]
    As \(\ind{G}{N} < \infty\), both kernels have finite index in \(H\), hence so does \(C\). Hence
    \[ \ind{G}{K} = \ind{AH\varphi}{BC\varphi} \leq \ind{A}{B} \ind{H\varphi}{C\varphi} \leq \ind{A}{B}\ind{H}{C}< \infty,\]
    so \(G\) is nilpotent-by-finite as well.
\end{proof}

With this final \zcref*[noref,nocap]{lem:ANquint3} proven, we are ready to construct the required algorithm. Since we already proved the existence of an algorithm that solves \zcref{prob:mainproblem} when \(G\) is nilpotent-by-finite in the previous section, this takes very little work.

\begin{proof}[Proof of \zcref{alg:metabelian}]
    We first calculate whether or not \(G\) is nilpotent-by-finite, by constructing the Fitting subgroup of \(G\) and then checking whether its index in \(G\) is finite or not. If it is not, then by the contrapositive of \zcref{lem:ANquint3} \(R_A(\varphi,\psi) = \infty\). If it is, then we defer to \zcref{alg:nilpotentbyfinite}.
\end{proof}

\section{Abelian subgroup commuting with the derived subgroup}
\label{sec:abcd}
Once again, the algorithm below is the main result of this section. The additional conditions we place on the standard quintuple here are essentially the same conditions Roman'kov used in his proofs of \zcref{alg:tcp,alg:coin}, see \cite{roma10-a,roma21-a}.

\begin{algothm}
    \label{alg:commutingderived}
    There exists an algorithm that solves \zcref{prob:mainproblem} for every standard quintuple \((G,H,\varphi,\psi,A)\) where
    \begin{enumerate}
    \item \(A\) is abelian,
    \item \([A,G'] = 1\).
\end{enumerate}
    Such a quintuple will be called an \emph{ABCD-quintuple}.
\end{algothm}

To begin, we impose two additional conditions on the quintuple, which allow us to reduce the problem to that of a metabelian quintuple and hence apply the results from \zcref{sec:metabelian}.
\begin{algothm}
    \label{alg:ccstep1}
    There exists an algorithm that solves \zcref{prob:mainproblem} for every ABCD-quintuple \((G,H,\varphi,\psi,A)\) where \(G = AH\varphi\) and \(H' \leq \Coin(\varphi,\psi)\).
\end{algothm}
\begin{proof}   
    Let \(p\colon H \to \ab{H}\) be the natural projection to the abelianisation. Since \(A\) and \(G'\) commute, the map
    \[ A \times \ab{H} \to A \colon (a,\bar{h}) \mapsto a^{\bar{h}} \coloneq (h\varphi)^{-1}a(h\varphi),\]
    with \(h \in \bar{h}p^{-1}\), defines a well-defined action of \(\ab{H}\) on \(A\). Construct the semi-direct product \(S \coloneq \ab{H} \ltimes A\)
    and define the group homomorphisms
    \begin{align*}
        &\lambda \colon \ab{H} \to S \colon \bar{h} \mapsto (\bar{h},1),\\
        &\mu \colon \ab{H} \to S \colon \bar{h} \mapsto (\bar{h},(h\varphi)^{-1}(h\psi)),
    \end{align*}
    where again \(h \in \bar{h}p^{-1}\). Then \((S,\ab{H},\lambda,\mu,A)\) is a metabelian quintuple and the map
    \[ {\R}_A[\varphi,\psi] \to {\R}_A[\lambda,\mu]\colon [a]_{\varphi,\psi} \mapsto [a]_{\lambda,\mu}\]
    is a bijection, so it now suffices to apply \zcref{alg:metabelian}.
\end{proof}

Next, we eliminate the condition \(H' \leq \Coin(\varphi,\psi)\) from the above \zcref*[noref,nocap]{alg:ccstep1}.

\begin{algothm}
    \label{alg:ccstep2}
    There exists an algorithm that solves \zcref{prob:mainproblem} for every ABCD-quintuple \((G,H,\varphi,\psi,A)\) where \(G = AH\varphi\).
\end{algothm}
\begin{proof}
    Define the map
    \[ \delta \colon H' \to A\colon h \mapsto (h\varphi)^{-1}(h\psi),\]
    which is a group homomorphism since \(A\) and \(G'\) commute.
    It follows from the condition \(G = AH\varphi\) that the image \(H'\delta\) is a normal subgroup of \(G\).
    Let \(p \colon G \to G/H'\delta\) be the natural projection and set \(\bar{G} \coloneq Gp\), \(\bar{A} \coloneq Ap\), \(\bvarphi \coloneq \varphi p\) and \(\bpsi \coloneq \psi p\). Then the map
    \[{\R}_A[\varphi,\psi] \to {\R}_{\bar{A}}[\bvarphi,\bpsi] \colon [a]_{\varphi,\psi} \mapsto [ap]_{\bvarphi,\bpsi}\]
    is a bijection. Hence, it is sufficient to solve \zcref{prob:mainproblem} for \((\bar{G},H,\bvarphi,\bpsi,\bar{A})\), which is an ABCD-quintuple satisfying the requirements of \zcref{alg:ccstep1}.
\end{proof}

Finally, we eliminate the condition \(G = A H\varphi\) from the above \zcref*[noref,nocap]{alg:ccstep2}.

\begin{proof}[Proof of \zcref{alg:commutingderived}]
    Set \(K \coloneq AH\varphi\) and define
    \[ \lambda\colon H \to K\colon h \mapsto h\varphi, \quad \mu \colon H \to K\colon h \mapsto h\psi.\]
    Now \((K,H,\lambda,\mu,A)\) is an ABCD-quintuple with \(K = AH\lambda\). Then the map
    \[{\R}_A[\varphi,\psi] \to {\R}_{A}[\lambda,\mu] \colon [a]_{\varphi,\psi} \mapsto [a]_{\lambda,\mu}\]
    is bijective, so it suffices to apply \zcref{alg:ccstep2} to \((K,H,\lambda,\mu,A)\).
\end{proof}

\section{General case}
\label{sec:general}
With most of the preliminary algorithms now at our disposal, we are ready to tackle \zcref{prob:mainproblem} in full generality.

\begin{algothm}
    \label{alg:nilpotentbyabelian}
    There exists an algorithm that solves \zcref{prob:mainproblem} for every standard quintuple \((G,H,\varphi,\psi,N)\) where \(G\) is nilpotent-by-abelian.
\end{algothm}
\begin{proof}
    Since \(G\) is nilpotent-by-abelian, its derived subgroup \(G'\) is nilpotent, say of class \(c\). We proceed by induction on \(c\). If \(c = 0\), then \(G\) is abelian and \(N\) is central, hence we defer to \zcref{alg:central}.
    
    Now suppose that \(c > 0\). As in the proof of \zcref{alg:nilpotent}, we will apply \zcref{alg:decomp}, but for the subgroup \(K \coloneq Z(G')\). If we consider the quintuple \((\bar{G},H,\bvarphi,\bpsi,\bar{N})\), then \(\bar{G}\) is nilpotent-by-abelian with \(\bar{G}'\) nilpotent of class \(c-1\), so the required algorithm exists by induction. For every \(n \in N\), the quintuple \((G,C_n,\varphi_n,\psi_n,M)\) has \(M = K \cap N\) abelian and \([M,G'] \leq [K,G'] = 1\), hence it is an ABCD-quintuple and we can apply \zcref{alg:commutingderived}.
\end{proof}

\begin{proof}[Proof of \zcref{alg:mainalgA}]
    Let \(K\) be a nilpotent-by-abelian finite index normal subgroup of \(G\). \zcref{alg:nilpotentbyabelian} can solve \zcref{prob:mainproblem} for every standard quintuple of the form \((K,L,\lambda,\mu,M)\), hence the result follows from \zcref{alg:byfinite}.
\end{proof}

\section{Summary of the algorithms}
\label{sec:summary}

In this section, we condense the algorithms from the previous sections into a few short snippets of pseudocode. For the convenience of the reader, we will also indicate with comments where the algorithms from the preceding sections begin.

Specifically, we provide pseudocode for the functions \textsc{RepsNil}, \textsc{RepsNBF} and \textsc{RepsNBA}, which solve \zcref{prob:mainproblem} for nilpotent, nilpotent-by-finite and nilpotent-by-abelian groups, respectively. Each of these functions takes as input a standard quintuple \((G,H,\varphi,\psi,N)\), and returns either a finite set of representatives of \(\R_N[\varphi,\psi]\), or ``\(\infty\)'' if the Reidemeister number \(R_N(\varphi,\psi)\) is infinite. To keep the code clear and concise, we will omit the definitions of homomorphisms induced by the arguments \(\varphi\) and \(\psi\).

We start with \textsc{RepsNil}, which corresponds to \zcref{alg:nilpotent} and thus assumes that the argument \(G\) is nilpotent.
\begin{algorithm}
\begin{algorithmic}[1]
    \Function{RepsNil}{$G,H,\varphi,\psi,N$}
    \If{\(N = 1\)}
        \Return \(\{1\}\)
    \EndIf
    \State Set \(K \coloneq Z(G)\) \Comment{\zcref{alg:decomp}}
    \State Set \(p \colon G \to G/K\), \(\bar{N} \coloneq Np\), \(\bar{G} \coloneq Gp\) 
    \State Run \textsc{RepsNil}(\(\bar{G},H,\bvarphi,\bpsi,\bar{N}\)) \Comment{\zcref{alg:nilpotent}}
    \If{\(R_{\bar{N}}(\bvarphi,\bpsi) = \infty\)}
        \Return \(\infty\)
    \EndIf
    \State Set \(\bar{n}_1, \ldots, \bar{n}_r\) representatives of  \(\R_{\bar{N}}[\bvarphi,\bpsi]\)
    \State Compute preimages \(n_1, \ldots, n_r\)
    \State Set \(A \coloneq N \cap K\)
    \For{$i = 1, \ldots, r$}
        \State Set \(H_i \coloneq \Coin_K(\varphi\iota_{n_i},\psi)\) \Comment{\zcref{alg:coin}}
        \State Set \(\delta_i \colon H_i \to A \colon h \mapsto (h\varphi\iota_{n_i})^{-1}(h\psi)\) \Comment{\zcref{alg:central}}
        \State Set \(q_i\colon A \to A/H_i\delta_i\), \(\tilde{A}_i \coloneq Aq_i\) 
        \If{\(\tilde{A}_i\) is infinite}
            \Return \(\infty\)
        \EndIf
        \State Enumerate \(\tilde{A}_i\) by \(\tilde{a}_{i1}, \ldots, \tilde{a}_{is_i}\)
        \State Compute preimages \(a_{i1}, \ldots, a_{is_i}\)
    \EndFor
    \State \Return \( \set{ n_ia_{ij}}{ i = 1,\ldots,r; \  j = 1,\ldots, s_i}\)
    \EndFunction
\end{algorithmic}
\end{algorithm}

Next is the function \textsc{RepsNBF}, which requires its argument \(G\) to be nilpotent-by-finite and hence corresponds to \zcref{alg:nilpotentbyfinite}. Note that we assume the existence of a function \textsc{RepsFin} that solves \zcref{prob:mainproblem} in the case that \(G\) is finite. Such a function can be constructed using standard orbit-stabiliser algorithms (see e.g.\@ Sections 4.1 and 8.6 in \cite{heo05-a}).
\begin{algorithm}
\begin{algorithmic}[1]
    \Function{RepsNBF}{$G,H,\varphi,\psi,N$}
    \State Set \(K \coloneq \Fitt(G)\), \(L \coloneq K\varphi^{-1} \cap K\psi^{-1}\) \Comment{\zcref{alg:byfinite}}
    \State Set \(p \colon G \to G/K\), \(\bar{N} \coloneq Np\), \(\bar{G} \coloneq Gp\)
    \State Run \textsc{RepsFin}(\(\bar{G},H,\bvarphi,\bpsi,\bar{N}\))
    \State Set \(\bar{n}_1, \ldots, \bar{n}_r\) representatives of  \(\R_{\bar{N}}[\bvarphi,\bpsi]\)
    \State Compute preimages \(n_1, \ldots, n_r\)
    \State Set \(A \coloneq N \cap K\)
    \For{$i = 1, \ldots, r$}
        \State Set \(L_i \coloneq \Coin_K(\restr{\varphi}{L}\iota_{n_i},\restr{\psi}{L})\) \Comment{\zcref{alg:coin}}
        \State Run \textsc{RepsNil}(\(K,L_i,\varphi_i,\psi_i,A\)) \Comment{\zcref{alg:nilpotent}}
        \If{\(R_A(\varphi_i,\psi_i) = \infty\)}
            \Return \(\infty\)
        \EndIf
        \State Set \(a_{i1}, \ldots, a_{is_i}\) representatives of \(\R_{A}[\varphi_i,\psi_i]\)
        \State Reduce to representatives \(a'_{i1}, \ldots, a'_{it_i}\)  of  \(\R_{A}[\varphi\iota_{n_i},\psi]\) \Comment{\zcref{alg:tcp}}
    \EndFor
    \State \Return \( \set{ n_ia'_{ij}}{ i = 1,\ldots,r; \  j = 1,\ldots, t_i}\)
    \EndFunction
\end{algorithmic}
\end{algorithm}

Finally, we have the function \textsc{RepsNBA}, which requires its argument \(G\) to be nilpotent-by-abelian and corresponds to \zcref{alg:nilpotentbyabelian}.

\begin{algorithm}
\begin{algorithmic}[1]
    \Function{RepsNBA}{$G,H,\varphi,\psi,N$}
    \If{\(N = 1\)}
        \Return \(\{1\}\)
    \EndIf
    \State Set \(K \coloneq Z(G')\) \Comment{\zcref{alg:decomp}}
    \State Set \(p \colon G \to G/K\), \(\bar{N} \coloneq Np\), \(\bar{G} \coloneq Gp\)
    \State Run \textsc{RepsNBA}(\(\bar{G},H,\bvarphi,\bpsi,\bar{N}\)) \Comment{\zcref{alg:nilpotentbyabelian}}
    \If{\(R_{\bar{N}}(\bvarphi,\bpsi) = \infty\)}
        \Return \(\infty\)
    \EndIf
    \State Set \(\bar{n}_1, \ldots, \bar{n}_r\) representatives of  \(\R_{\bar{N}}[\bvarphi,\bpsi]\)
    \State Compute preimages \(n_1, \ldots, n_r\)
    \State Set \(A \coloneq N \cap K\)
    \For{$i = 1, \ldots, r$}
        \State Set \(H_i \coloneq \Coin_K(\varphi\iota_{n_i},\psi)\), \(G_i \coloneq AH_i\varphi\iota_{n_i}\) \Comment{\zcref[nosort]{alg:coin,alg:commutingderived}}
        \State Set \(\delta_i \colon H_i' \to A \colon h \mapsto (h\varphi\iota_{n_i})^{-1}(h\psi)\) \Comment{\zcref{alg:ccstep2}}
        \State Set \(q_i\colon G_i \to G_i/H_i'\delta_i\), \(\tilde{A}_i \coloneq Aq_i\), \(\tilde{G}_i \coloneq \ab{H}_i \ltimes \tilde{A}_i\)  \Comment{\zcref{alg:ccstep1}}
        \If{\(\tilde{G}_i\) is not nilpotent-by-finite}
            \Return \(\infty\) \Comment{\zcref{alg:metabelian}}
        \EndIf
        \State Run \textsc{RepsNBF}(\(\tilde{G}_i,\ab{H}_i,\tvarphi_i,\tpsi_i,\tilde{A}_i\)) \Comment{\zcref{alg:nilpotentbyfinite}}
        \If{\(R_{\tilde{A}_i}(\tvarphi_i,\tpsi_i) = \infty\)}
            \Return \(\infty\)
        \EndIf
        \State Set \(\tilde{a}_{i1}, \ldots, \tilde{a}_{is_i}\) representatives of  \(\R_{\tilde{A}_i}[\tilde{\varphi}_i,\tilde{\psi}_i]\)
        \State Compute preimages \(a_{i1}, \ldots, a_{is_i}\)
    \EndFor
    \State \Return \( \set{ n_ia_{ij}}{ i = 1,\ldots,r; \  j = 1,\ldots, s_i}\)
    \EndFunction
\end{algorithmic}
\end{algorithm}

We have omitted pseudocode for \zcref{alg:mainalgA}, since this would be identical to \textsc{RepsNBF} with only two modifications: in line \(2\) we take \(K\) to be a finite index normal nilpotent-by-abelian subgroup of \(G\) instead of \(\Fitt(G)\), and in line \(10\) we run \textsc{RepsNBA} instead of \textsc{RepsNil}.

To finish this section, we briefly remark that \zcref{alg:mainalgA} has one major (potential) bottleneck: having to solve the twisted conjugacy problem many times. This happens every time we use \zcref{alg:byfinite}, in particular when calculating the image of the map \(\pi\) described in the proof of said \zcref[noref,nocap]{alg:byfinite}. Note that this bottleneck is avoided when the argument \(G\) is nilpotent.

\section{Double cosets}
\label{sec:doublecosets}
If \(G\) is a group with subgroups \(U, V\) and element \(g\), then the \emph{double coset} \(UgV\) is defined as \(\set{ ugv }{ u \in U, v \in V}\). The set of all \((U,V)\)-double cosets is usually denoted by \(U\backslash G / V\), and the number of \((U,V)\)-double cosets is called the \emph{double coset index} of the pair \((U,V)\).

As for twisted conjugacy classes, there is both a membership and a search problem for double cosets.

\begin{problem}
    \label{prob:dcmp}
    Given a group \(G\), elements \(x\) and \(y\) of \(G\) and \(U,V \leq G\), determine whether \(x \in UyV\).
\end{problem}

\begin{problem}
    \label{prob:dcsp}
    Given a group \(G\), elements \(x\) and \(y\) of \(G\) and \(U,V \leq G\) such that \(x \in UyV\), find \(u \in U\) and \(v \in V\) such that \(x = uyv\).
\end{problem}

In addition to studying a single double coset, we are also interested in determining how a group is partitioned by its double cosets.

\begin{problem}
    \label{prob:dcip}
    Given a group \(G\) and \(U,V \leq G\), compute the double coset index of the pair \((U,V)\).
\end{problem}

\begin{problem}
    \label{prob:dcrp}
    Given a group \(G\) and \(U,V \leq G\) with finite double coset index, find \(g_1, \ldots, g_r \in G\) such that \(G =  Ug_1V \sqcup Ug_2V \sqcup \cdots \sqcup Ug_rV\).
\end{problem}

\begin{algothm}
    There exist algorithms that solve \zcref{prob:dcmp,prob:dcsp,prob:dcip,prob:dcrp} for polycyclic-by-finite groups.
\end{algothm}
\begin{proof}
Set \(H \coloneq U \times V\) and consider the homomorphisms
\[
\varphi \colon H \to G \colon (u,v) \mapsto u, \qquad \psi \colon H \to G \colon (u,v) \mapsto v,
\]
then \([g]_{\varphi,\psi} = UgV\) for every \(g \in G\). Therefore, the problems in question can be solved using \zcref{alg:mainalgA,alg:tcp}.
\end{proof}

\section{Affine actions}
\label{sec:affineactions}
In \zcref{sec:metabelian} we introduced the notions of group modules and derivations. While commutativity is a natural condition to impose in the context of group cohomology, we can extend the definition of a derivation to having non-abelian range.

\begin{definition}
Let \(G,H\) be groups and let \(H\) act on \(G\) via automorphisms. A map \( \delta \colon H \to G \) is called a \emph{derivation} if \((h_1h_2)\delta = (h_1\delta)^{h_2}(h_2\delta)\) for all \(h_1,h_2 \in H\).
\end{definition}

To a derivation we can associate a new action: the \emph{affine action} associated to the derivation.

\begin{definition}
    Let a group \(H\) act on a group \(G\) via \(\lambda \colon H \to \Aut(G) \colon h \mapsto \lambda_h\), and let \(\delta \colon H \to G\) be a derivation with respect to this action. The \emph{affine action} associated to \(\delta\) is the action given by
    \[ \alpha \colon G \times H \to G \colon (g,h) \mapsto (g\lambda_h)(h\delta).\]
\end{definition}

We can see this as a natural extension of ``acting via automorphisms''. Rather than a map from \(H\) to \(\Aut(G)\), this action now corresponds to a map from \(H\) to \(\Aff(G) \coloneq \Aut(G) \ltimes G\), i.e.\@ the map
\[H \to \Aff(G) \colon h \mapsto (\lambda_h,h\delta).\]
An action of \(H\) on \(G\) induced by a homomorphism \(H \to \Aff(G)\) will also be called an affine action.

\begin{proposition}
    \label{prop:twiconactisaffine}
    Let \(G\) and \(H\) be groups and let \(\varphi,\psi \in \Hom(H,G)\). Then the \((\varphi,\psi)\)-twisted conjugation action of \(H\) on \(G\) is an affine action.
\end{proposition}
\begin{proof}
    The affine action induced by the group homomorphism
    \[
    H \to \Aff(G) \colon h \mapsto (\iota_{h\varphi},(h\varphi)^{-1}(h\psi))
    \]
    coincides with the \((\varphi,\psi)\)-twisted conjugation action.
\end{proof}

The converse of \zcref{prop:twiconactisaffine} is not true: in general, an affine action of a group \(H\) on a group \(G\) need not coincide with the twisted conjugation action of some pair of homomorphisms.
However, using a construction very similar to that of \zcref{alg:ccstep1}, we can still leverage algorithms designed for twisted conjugacy to do calculations for affine actions.

\begin{algothm}
    \label{alg:affineaction}
    Let \(G\) and \(H\) be polycyclic-by-finite groups and let \(\alpha \colon H \to \Aff(G)\) be a group homomorphism. There exist algorithms that, for the affine action \((g,h) \mapsto g^h\) induced by \(\alpha\):
    \begin{enumerate}
        \item calculate the stabiliser of an element \(g\);
        \item decide whether two elements \(g_1,g_2\) belong to the same orbit;
        \item calculate an element \(h\) such that \({g_1}^h = g_2\), if it is known to exist;
        \item calculate the number of orbits;
        \item find representatives of the orbits if there are only finitely many.
    \end{enumerate}
\end{algothm}
\begin{proof}
    Suppose that the homomorphism \(\alpha\) is given by
    \[ \alpha \colon H \to \Aut(G) \ltimes G \colon h \mapsto (\lambda_h,g_h).\]
    Then \(H\) also acts on \(G\) via \(\lambda \colon H \to \Aut(G) \colon h \mapsto \lambda_h\), and we can construct the
    semi-direct product \(K \coloneq H \ltimes_\lambda G\) which is polycyclic-by-finite. Consider the homomorphisms
    \begin{align*}
        \varphi &\colon H \to K \colon h \mapsto (h,1_G),\\
        \psi &\colon H \to K \colon h \mapsto (h,g_h).
    \end{align*}
    Identifying \(G\) with its inclusion in \(K\), we note that \(g^h = (h\varphi)^{-1}g(h\psi)\) for every \(g \in G\) and every \(h \in H\). This means that the orbit and stabiliser of \(g\) under the affine action coincide with the orbit and stabiliser of \(g\) under the \((\varphi,\psi)\)-twisted conjugation action, respectively. Thus, we can apply \zcref{alg:mainalgA,alg:tcp,alg:coin}, noting that the stabiliser of \(g\) is exactly the coincidence group \(\Coin(\varphi\iota_g,\psi)\).
\end{proof}

As we originally introduced affine actions using derivations, it should come as no surprise that we can use these algorithms to study derivations as well.

\begin{algothm}
    \label{alg:derivations}
    Let \(G\) and \(H\) be polycyclic-by-finite groups and let \(\delta\colon H \to G\) be a derivation. There exist algorithms that:
    \begin{enumerate}[ref=(\arabic*)]
        \item\label{alg:derivations:ker} calculate the kernel of \(\delta\);
        \item\label{alg:derivations:img} decide whether an element \(g \in G\) lies in the image of \(\delta\);
        \item\label{alg:derivations:preimg} calculate a preimage \(h \in g\delta^{-1}\) for \(g \in \im(\delta)\);
        \item\label{alg:derivations:fullpreimg} calculate the full set of preimages \(g\delta^{-1}\) of \(g \in G\);
        \item\label{alg:derivations:surj} decide whether \(\delta\) is surjective.
    \end{enumerate}
\end{algothm}
\begin{proof}
    If we consider the affine action associated to \(\delta\), then \(\im(\delta)\) and \(\ker(\delta)\) are the orbit and stabiliser of the identity \(1_G\), respectively. Thus the existence of algorithms \zcref[noname]{alg:derivations:ker,alg:derivations:img,alg:derivations:preimg} follows immediately from \zcref{alg:affineaction}. We discuss the remaining two algorithms.
    \begin{enumerate}[start=4]
        \item If \(g \notin \im(\delta)\), which we can verify using \zcref[noname]{alg:derivations:img}, then the set of preimages is empty. If \(g \in \im(\delta)\), then we can calculate \(K \coloneq \ker(\delta)\) using \zcref[noname]{alg:derivations:ker} and find some \(h \in g\delta^{-1}\) using \zcref[noname]{alg:derivations:preimg}. The full set of preimages \(g\delta^{-1}\) is exactly the right coset \(Kh\).
        \item Observe that \(\delta\) is surjective if and only if \([1]_{\varphi,\psi} = G\), and this is equivalent to \(R_G(\varphi,\psi)\) being \(1\). Therefore it suffices to compute \(R_G(\varphi,\psi)\) using \zcref{alg:mainalgA}. \qedhere
    \end{enumerate}
\end{proof}

We remark that \zcref*{alg:derivations}\zcref[noname]{alg:derivations:surj} generalises \cite[Thm.~B]{robi02-a}, which states that there exists an algorithm that can decide whether a derivation from a polycyclic-by-finite group to a finitely generated abelian group is surjective. \zcref[cap]{alg:affineaction,alg:derivations} have been implemented in the \textsf{GAP} package \texttt{TwistedConjugacy} \cite{tert26-a}.

\section*{Acknowledgements}
The author would like to thank Karel Dekimpe for the helpful discussions, Marlies Vantomme for her useful remarks and suggestions, and the anonymous referees for their careful reading and their valuable comments, corrections and suggestions.
{\emergencystretch=1em
\printbibliography}

\end{document}